\documentstyle{amsppt}
\magnification=1200
\TagsOnRight
\NoBlackBoxes
\topmatter
\NoBlackBoxes
\title Finitely Smooth Reinhardt Domains\\ 
with Non-Compact\\ 
Automorphism Group
\endtitle
\rightheadtext{Finitely Smooth Reinhardt Domains}
\footnote[]{{\bf Mathematics
Subject Classification:} 32A07, 32H05, 32M05 \hfill}
\footnote[]{{\bf Keywords and Phrases:} Automorphism
groups, Reinhardt domains. \hfill}
\author A. V. Isaev \ \ \ and \ \ \ S. G. Krantz
\endauthor 
\abstract {We give a complete description of bounded
Reinhardt domains of finite boundary smoothness that have
non-compact automorphism group. As part of this
program, we show that the classification of domains with
non-compact automorphism group and having only finite
boundary smoothness is considerably more complicated than
the classification of such domains that have infinitely
smooth boundary.}
\endabstract
     
\endtopmatter   
\document
\def\qed{{\hfill{\vrule height7pt width7pt
depth0pt}\par\bigskip}}

Let $D\subset{\Bbb C}^n$ be a bounded domain, and suppose that the
group $\text{Aut}(D)$ of holomorphic automorphisms of $D$ is
non-compact in the topology of uniform convergence on compact subsets
of $D$. This means that there exist points $q\in\partial D$, $p\in D$ and a
sequence $\{f_j\}\subset \text{Aut}(D)$ such that $f_j(p)\rightarrow
q$ as $j\rightarrow \infty$. 

We also assume that $D$ is a Reinhardt domain, i.e.\ that the
standard action of the $n$-dimensional torus ${\Bbb T}^n$ on ${\Bbb
C}^n$,
$$
z_j\mapsto e^{{i\phi}_j}z_j,\qquad {\phi}_j\in {\Bbb R},\quad
j=1,\dots,n,
$$
leaves $D$ invariant.

In \cite{FIK1} we gave a complete classification of bounded Reinhardt 
domains with 
non-compact automorphism group and $C^\infty$-smooth
boundary. For the sake of completeness we quote the main
result of \cite{FIK1} below:

\proclaim{Theorem 1} If $D$ is a bounded Reinhardt domain
in ${\Bbb C}^n$ with $C^{\infty}$-smooth boundary, and if 
$\text{\rm Aut}(D)$
is not compact then, up to dilations and permutations of coordinates, 
$D$ is  
a domain of the form
$$
\left\{|z^1|^2+\sum_{j=2}^p|z^j|^{2m_j}+P(|z^2|,\dots,|z^p|)<1\right\},
$$
where $P$ is a polynomial:
$$
P(|z^2|,\dots,|z^p|)=\sum_{l_2,\dots,l_p}a_{l_2,\dots,l_p}
|z^2|^{2l_2}\dots
|z^p|^{2l^p},\tag{1}
$$
$a_{l_2,\dots,l_p}$ are real parameters, $m_j\in{\Bbb N}$, with the
sum taken over all $(p-1)$-tuples $(l_2,\dots,l_p)$, $l_j\in{\Bbb Z}$,
$l_j\ge 0$, where at least two entries are non-zero, such that
$\sum_{j=2}^p\frac{l_j}{m_j}=1$, and the complex variables
$z_1,\dots,z_n$ are divided into $p$ non-empty groups $z^1,\dots,z^p$. 
In addition, the polynomial 
$$
\tilde P(|z^2|,\dots,|z^p|)=\sum_{j=2}^p|z^j|^{2m_j}+P(|z^2|,\dots,|z^p|)
$$
is non-negative in ${\Bbb C}^{n-n_1}$ ($n_1$ is the number of variables 
in the group $z^1$), and the domain
$$
\left\{(z^2,\dots,z^p)\in{\Bbb C}^{n-n_1}\: 
\tilde P(|z^2|,\dots,|z^p|)<1\right\}
$$
is bounded.
\endproclaim

In this paper we generalize Theorem 1 to the case when the
boundary of the domain is only $C^k$-smooth, $k\ge 1$. To the best of 
our knowledge, this is the first attempt to obtain a general result for
bounded domains with non-compact automorphism group and
boundary of finite smoothness. 

First of all, we note that, up to a certain point, the proof of 
Theorem 1 in \cite{FIK1} is valid for domains with only $C^1$-smooth
boundary. The $C^\infty$-assumption was only used in Lemmas 1.6
and 1.8 of \cite{FIK1}. Therefore, the proof of Theorem 1 in
\cite{FIK1} also gives the following proposition:

\proclaim{Proposition 2} If $D$ is a bounded Reinhardt domain in
${\Bbb C}^n$ with $C^k$-smooth boundary, $k\ge 1$, and if $\text{\rm
Aut}(D)$ is not compact, then, by suitable dilations and permutations of
coordinates, the domain $D$ is equivalent to a Reinhardt domain $G$ such 
that:
\smallskip

\noindent (i) The set $A$ of all points $(z^1,0,\dots,0)$, with $|z^1|=1$,
lies in $\partial G$.
\smallskip

\noindent (ii) In a neighbourhood of $A$, $G$ is written in the form 
$$ 
\left\{|z^1|^2+\phi(|z^2|,\dots,|z^p|)<1\right\},\tag{2} 
$$ 
where $\phi(x_2,\dots,x_p)$ is a non-negative $C^k$-smooth function in a
neighbourhood of the origin in ${\Bbb R}^{p-1}$ such that
$\phi(|z^2|,\dots,|z^p|)$ is also $C^k$-smooth in a neighbourhood of the
origin in ${\Bbb C}^{n-n_1}$, and such that 
$$
\phi\left(t^{\frac{1}{\alpha_2}}x_2,\dots,
t^{\frac{1}{\alpha_p}}x_p\right)=t\phi(x_2,\dots,x_p)\tag{3}
$$ 
near the origin in ${\Bbb R}^{p-1}$ for $1\le t\le 1+\epsilon$ and
some $\epsilon>0$. Here $\alpha_j>0$, $j=2,\dots,p$, and each
$\alpha_j$ is either an even integer or, if it is not an even integer,
then $\alpha_j>2k$. In addition, the function $\phi$ satisfies 
$$
\phi(0,\dots,0,|z^j|,0,\dots,0)=|z^j|^{\alpha_j},
$$
for $j=2,\dots,p$.
\smallskip

\noindent (iii) $G$ has the form
$$
G=\left\{(z^1,\dots,z^p)\in{\Bbb C}^n\:|z^1|<1,
\left(\frac{z^2}{(1-|z^1|^2)^{\frac{1}{\alpha^2}}},\dots,
\frac{z^p}{(1-|z^1|^2)^{\frac{1}{\alpha^p}}}\right)\in\tilde
G\right\},\tag{4}
$$
where $\tilde G$ is a bounded Reinhardt domain in ${\Bbb C}^{n-n_1}$.
\endproclaim
\smallskip

We are now going to derive from Proposition 2 the 
following theorem, which is the main result of the present
paper.

\proclaim{Theorem 3} If $D$ is a bounded Reinhardt domain in
${\Bbb C}^n$ with $C^k$-smooth boundary, $k\ge 1$, and if $\text{\rm
Aut}(D)$ is not compact, then, up to dilations and permutations of
coordinates, $D$ is a domain of the form
$$ 
\left\{|z^1|^2+\psi(|z^2|,\dots,|z^p|)<1\right\},\tag{5} 
$$ 
where $\psi(x_2,\dots,x_p)$ is a non-negative $C^k$-smooth function in 
${\Bbb R}^{p-1}$ such that \break
$\psi(|z^2|,\dots,|z^p|)$ is $C^k$-smooth
in ${\Bbb C}^{n-n_1}$, and  
$$
\psi\left(t^{\frac{1}{\alpha_2}}x_2,\dots,
t^{\frac{1}{\alpha_p}}x_p\right)=t\psi(x_2,\dots,x_p)\tag{6}
$$ 
in ${\Bbb R}^{p-1}$ for all $t\ge 0$. Here $\alpha_j>0$, $j=2,\dots,p$, 
and each
$\alpha_j$ is either an even integer or, if it is not an even integer,
then $\alpha_j>2k$. In addition, the function $\psi$ satisfies 
$$
\psi(0,\dots,0,|z^j|,0,\dots,0)=|z^j|^{\alpha_j},\tag{7}
$$
for $j=2,\dots,p$, and the domain
$$
\left\{(z^2,\dots,z^p)\in{\Bbb C}^{n-n_1}\: \psi(|z^2|,\dots,|z^p|)<1
\right\}\tag{8}
$$
is bounded.

 \endproclaim

\demo{Proof} First of all, using the weighted homogeneity 
property (3), we extend the function $\phi$ from a
neighbourhood of the origin (see (2)) to a $C^k$-smooth
function $\psi$ on ${\Bbb R}^{p-1}$. Consider the surface
$$
S_{\delta}=\left\{|x_2|^{\alpha_2}+\dots+|x_p|^{\alpha_p}=
\delta\right\},
$$ 
and choose $\delta>0$ such that $S_{\delta}$ lies in the
neighbourhood of the origin in ${\Bbb R}^{p-1}$ where
$\phi$ is defined and of class $C^k$ and where (3) holds
for $1\le t\le 1+\epsilon$. Further let 
$$
S_{\delta}^-=\left\{|x_2|^{\alpha_2}+\dots+
|x_p|^{\alpha_p}\le\delta\right\},
$$ and 
$$
S_{\delta}^+=\left\{|x_2|^{\alpha_2}+\dots+
|x_p|^{\alpha_p}\ge\delta\right\}.
$$

We now define the extension $\psi$ 
for any $x=(x_2,\dots,x_p)\in{\Bbb R}^{p-1}$ as follows:
$$
\psi(x)=\cases
\phi(x),\text{ if $x\in S_{\delta}^-$};\\
\frac{|x_2|^{\alpha_2}+\dots+|x_p|^{\alpha_p}}{\delta}\times\\
\times\phi\left(\frac{x_2\delta^{\frac{1}{\alpha_2}}}
{(|x_2|^{\alpha_2}+\dots+|x_p|^{\alpha_p})^{\frac{1}{\alpha_2}}},\dots,
\frac{x_p\delta^{\frac{1}{\alpha_p}}}{(|x_2|^{\alpha_2}+
\dots+|x_p|^{\alpha_p})^{\frac{1}{\alpha_p}}}\right),\text{
if $x\in S_{\delta}^+$}. \endcases\tag{9}
$$
Because of property (3), this definition implies that,
for any $x\ne 0$ 
$$
\multline
\psi(x)=\frac{|x_2|^{\alpha_2}+\dots+|x_p|^{\alpha_p}}{\delta}\times\\
\times\phi\left(\frac{x_2\delta^{\frac{1}{\alpha_2}}}{(|x_2|^{\alpha_2}+
\dots+|x_p|^{\alpha_p})^{\frac{1}{\alpha_2}}},\dots,
\frac{x_p\delta^{\frac{1}{\alpha_p}}}{(|x_2|^{\alpha_2}+
\dots+|x_p|^{\alpha_p})^{\frac{1}{\alpha_p}}}\right),
\endmultline\tag{10}
$$
and $\psi(0)=0$. Further, since, for $j=2,\dots,p$, $\alpha_j$ is 
either a positive integer, or, if not, $\alpha_j>2k$, one has that
$\psi\in C^k\left({\Bbb R}^{p-1}\right)$ and $\psi(|z^2|,\dots,|z^p|)\in
C^k\left({\Bbb C}^{n-n_1}\right)$. Next, (10) implies that $\psi$
has property (6) for all $x\in {\Bbb R}^{p-1}$ and $t\ge
0$, as  well
as property (7). It is also clear that $\psi\ge 0$.

We will now show that the domain $G$ has the 
form (5), with $\psi$ defined in (9). Let $U$ be a
neighbourhood of the set $A$ (see (i) of Proposition 2)
such that $G\cap U$ is given by (2). We can assume that
$U=V\times W$, where $V$ is a neighbourhood of the unit
sphere in ${\Bbb C}^{n_1}$, and $W$ is a neighbourhood of
the origin in ${\Bbb C}^{n-n_1}$. Take $\sigma>0$ and
consider $G_{\sigma}=G\cap\{|z^1|^2=1-\sigma\}$. Since
$\tilde G$ is bounded, representation (4) implies that, if
$\sigma$ is sufficiently small, $G_{\sigma}\subset U$, and
$\overline{G_{\sigma}}$ is a compact subset of $U$. It
then follows from (2) that $G_{\sigma}$ is given by $$
G_{\sigma}=\left\{(z^1,\dots,z^p)\in
U\:|z^1|^2=1-\sigma,\,\,
\phi(|z^2|,\dots,|z^p|)<\sigma\right\}, $$ and the set
$\{(z^2,\dots,z^p)\in W\:\phi(|z^2|,\dots,|z^p|)\le
\sigma\}$ is compact in $W$.

Further, since the extension $\psi$ of $\phi$ has property (6), 
$G_{\sigma}$ can be rewritten as
$$
G_{\sigma}=\left\{(z^1,\dots,z^p)\in {\Bbb C}^n\:|z^1|^2=1-\sigma,\,\,
\psi(|z^2|,\dots,|z^p|)<\sigma\right\}.
$$
On the other hand, (4) gives
$$
G_{\sigma}=\left\{(z^1,\dots,z^p)\in{\Bbb C}^n\:|z^1|^2=1-\sigma,\,\,
\left(\frac{z^2}{\sigma^{\frac{1}{\alpha^2}}},\dots,
\frac{z^p}{\sigma^{\frac{1}{\alpha^p}}}\right)\in\tilde
G\right\},
$$
which implies that
$$
\tilde G=\left\{(z^2,\dots,z^p)\in {\Bbb C}^{n-n_1}\:
\psi(\sigma^{\frac{1}{\alpha_2}}|z^2|,\dots,\sigma^{\frac{1}{\alpha_p}}
|z^p|)<\sigma\right\}.
$$
It now follows from homogeneity property (6) for $\psi$ that
$$
\tilde G=\left\{(z^2,\dots,z^p)\in {\Bbb C}^{n-n_1}\:
\psi(|z^2|,\dots,|z^p|)<1\right\}.
$$
Now (4) and (6) imply that $G$ is in fact given by formula (5).

Finally, domain (8) is bounded since it coincides with $\tilde G$.

The theorem is proved.\qed
\enddemo

For Reinhardt domains in ${\Bbb C}^2$, one has 
either $p=1$ or $p=2$. If $p=1$, then domain (5) is the
unit ball. If $p=2$, then, because of (7), the function
$\psi$ from (5) has the form $\psi=|z_2|^{\alpha}$,
$\alpha>0$. This observation gives the following corollary.

\proclaim{Corollary 4} If $D$ is a bounded Reinhardt domain
in ${\Bbb C}^2$ with $C^k$-smooth boundary, $k\ge 1$, and if 
$\text{\rm Aut}(D)$
is not compact, then, up to dilations and permutations of coordinates, 
$D$ has the form
$$
\{|z_1|^2+|z_2|^{\alpha}<1\},\tag{11}
$$
where $\alpha>0$ and either is an even integer or, if it is not
an even  integer, then $\alpha>2k$.
\endproclaim

\noindent{\bf Remark.} Note that Corollary 4 is 
reminiscent of a result of Bedford/Pinchuk (see
\cite{BP1}): a pseudoconvex smoothly bounded domain in
${\Bbb C}^2$ with non-compact automorphism group and
boundary of finite type in the sense of Kohn must be
biholomorphic to a domain of the form (11) where {\it
$\alpha$ is an even integer}.  The results of
Bedford/Pinchuk, and related conjectures, are discussed in
more details at the end of this paper.
\medskip

Theorem 3 reduces the classification problem for Reinhardt domains with 
non-compact automorphism group and $C^k$-smooth boundary to the problem
of describing $C^k$-smooth functions $\psi$ as in (5) that satisfy
weighted homogeneity condition (6). For $p\ge 3$, one can
construct examples of such functions in the following manner.
Consider the following set of $(p-1)$-tuples $s=(s_2,\dots,s_p)$
$$ \aligned M=&\biggl\{s=(s_2,\dots,s_p)\in{\Bbb
R}^{p-1}\:s_j\ge 0;\,\,\text{each $s_j$ is either an even
integer, or,}\\  &\text{if it is not an even integer, then
$s_j>2k$; $s$ has at least two non-zero}\\  &\text{entries; and
} \sum_{j=2}^p\frac{s_j}{\alpha_j}=1\biggr\}.
\endaligned\tag{12} $$ Let $\mu$ be an arbitrary finite measure
on the set $M$. Then the function $$
\psi(|z^2|,\dots,|z^p|)=\sum_{j=2}^p|z^j|^{\alpha_j}+\int_{M}|z^2|^{s_2}
\dots |z^p|^{s_p}\, d\mu\tag{13}  $$ 
has all the properties as stated in Theorem
3 above, provided $\psi\ge 0$ and the corresponding domain (8) is bounded.

We now give an explicit non-trivial example of a function of the form (13).
\medskip

{\bf Example 5.} Consider the case of ${\Bbb C}^3$ and 
let $p=3$, i.e. $z^j=z_j$, $j=1,2,3$. Let $k=2$,
$\alpha_2=\alpha_3=9$. Then it follows from (12) that 
$$
M=\{(s_2,s_3)\in{\Bbb R}^2\: s_2=9-s_3,\,4\le s_3\le 5\}\cup\{(2,7)\}
\cup\{(7,2)\}.
$$
We interpret $M$ as a subset of ${\Bbb R}$ 
parametrized by $s_3$ and let $d\mu=ds$ be the usual
Lebesgue measure on ${\Bbb R}$. Then the function defined
by (13) becomes 
$$ 
\multline
\psi(|z_2|,|z_3|)=|z_2|^9+|z_3|^9+|z_2|^9\int_4^5 
\frac{|z_3|^s}{|z_2|^s}\,
ds=\\|z_2|^9+|z_3|^9+\frac{1}{\log|z_3|^2-\log|z_2|^2}
\left(|z_2|^4|z_3|^5-|z_2|^5|z_3|^4\right).
\endmultline\tag{14} 
$$
The last term in function (14) and its first and second
derivatives are  defined to be equal to zero whenever $z_2=0$,
or $z_3=0$, or $|z_2|=|z_3|$. 

One can check directly that function (14) is 
indeed non-negative, $C^2$-smooth, has an appropriate
homogeneity property (6) with $\alpha_2=\alpha_3=9$, and
the corresponding domain (8) is bounded. The Reinhardt
domain $D\subset {\Bbb C}^3$ given by 
$$
D=\left\{|z_1|^2+|z_2|^9+|z_3|^9+\frac{1}{\log|z_3|^2-\log|z_2|^2}
\left(|z_2|^4|z_3|^5-|z_2|^5|z_3|^4\right)<1\right\}
$$ 
is a bounded domain with non-compact automorphism group
and $C^2$-smooth boundary.

Similar examples can be constructed in any complex 
dimension for any $p\ge 3$ and $k\ge 1$. Note that there
is considerable freedom in choosing a measure 
$\mu$ in (13). \qed 
\medskip

It is a reasonable question whether any 
function $\psi$ as in Theorem 3 is given by formula (13)
for an appropriate choice of $\mu$. Note that, as shown in
\cite{FIK1}, this holds if $k=\infty$, in which case the
entries of $(p-1)$-tuples $s$ from the set $M$ can only
be even integers and thus formula (13) turns into a
polynomial (see (1)). However, as demonstrated by the
following example, in the case of finite smoothness one
can find functions that have the weighted homogeneity
property, but that are not given by integration
against a measure as in (13). 
\medskip

{\bf Example 6.} As in Example 5, let again $n=3$, $p=3$ 
and $k=2$. We set $\alpha_2=\alpha_3=8$. Then it follows
from (12) that 
$$
M=\{(2,6)\}\cup\{(4,4)\}\cup\{(6,2)\}.
$$
Since $M$ is finite, the integral in (13) turns into a 
finite sum, and all functions of the form (13) are
real-analytic. We are now going to present a $C^2$-smooth
function $\psi(|z_2|,|z_3|)$ that has property (6)
with $\alpha_2=\alpha_3=8$ and such that $\psi$ is not
necessarily real-analytic.

Let $g\in C^2({\Bbb R})$ be such that 
$g(0)=0$ and $g(x)=x^2$ for $|x|>1$. Then a direct
calculation shows that  $$
\psi(|z_2|,|z_3|)=|z_2|^8+|z_3|^8+|z_2|^8g
\left(\frac{|z_3|^2}{|z_2|^2}\right)
$$
is $C^2$-smooth (for the last term 
$|z_2|^8g\left(\frac{|z_3|^2}{|z_2|^2}\right)$ we set its
value and the values of its first and second derivatives
to be equal to zero whenever $z_2=0$). The above function
$\psi$ satisfies (6) with $\alpha_2=\alpha_3=8$, but it,
of course, is not real-analytic for any non-trivial choice
of $g$. Also, if $g\ge 0$, one has that $\psi\ge 0$, and
the corresponding domain (8) is bounded. The Reinhardt
domain $D\subset{\Bbb C}^3$, 
$$
D=\left\{|z_1|^2+|z_2|^8+|z_3|^8+|z_2|^8g\left(\frac{|z_3|^2}{|z_2|^2}
\right)<1\right\} \ ,
$$ 
is then also bounded and has a non-compact automorphism  group
and $C^2$-smooth boundary. Such an example can be given in
any complex dimension for any $p\ge 3$, $k\ge 1$. \qed
\medskip

Example 6 shows that, most probably, a nice 
description of finitely smooth functions with
weighted homogeneity property does not exist, at least in
the form of an explicit formula such as (13). Therefore,
Theorem 3 is likely to be the best possible
classification result that one can hope to obtain for
Reinhardt domains of finite smoothness.

It also may be noted that weighted homogeneous functions
may be constructed by specifying them on the set
$S_1=\left\{|x_2|^{\alpha_2}+\dots+|x_p|^{\alpha_p}=1\right\}$ and
then extending to all of the space by homogeneity as in the proof of
Theorem 3 above (see (10)).  Such a construction is useful in that
it reduces the smoothness question to {\bf (i)} checking
smoothness on $S_1$; and {\bf (ii)}
checking smoothness at the origin (smoothness elsewhere
is automatic).

Along the lines of the preceding
discussion, one can consider the following examples of
domains with non-compact automorphism group and
$C^k$-smooth boundary, $k\ge 1$, that are not necessarily
Reinhardt:   
$$ 
\left\{|z_1|^2+\psi(z_2,\dots,z_n)<1\right\},\tag{15}  
$$ 
where $\psi(z_2,\dots,z_n)$ is a $C^k$-smooth function in
${\Bbb C}^{n-1}$ and 
$$
\psi\left(t^{\frac{1}{\alpha_2}}z_2,\dots,
t^{\frac{1}{\alpha_n}}z_n\right)=|t|\psi(z_2,\dots,z_n)\tag{16}
$$  
in ${\Bbb C}^{n-1}$ for all $t\in{\Bbb C}$. Here
$\alpha_j>0$, $j=2,\dots,n$, and
$t^{\frac{1}{\alpha_j}}=e^{\frac{1}{\alpha_j}(\log|t|+i\text{arg}\,t)}$,
for $t\ne 0$,
where $-\pi<\text{arg}\,t\le \pi$. Also, to guarantee that
domain (15) is bounded, one can assume that $\psi\ge 0$
and the domain   
$$
\left\{(z_2,\dots,z_n)\in{\Bbb C}^{n-1}\:
\psi(z_2,\dots,z_n)<1\right\}  $$ is bounded
(cf. Theorem 3). 

For any domain $D$ of the form (15), $\text{Aut}(D)$ is indeed 
non-compact, since it contains the subgroup
$$
\aligned
&z_1\mapsto \frac{z_1-a}{1-\overline{a}z_1},\\
&z_j\mapsto \frac{(1-|a|^2)^{\frac{1}{\alpha_j}}z_j}
{(1-\overline{a}z_1)^{\frac{2}{\alpha_j}}},\qquad j=2,\dots,n,
\endaligned
$$
where $|a|<1$. In addition to the above automorphisms, domains (15) 
are also invariant under the special rotations
$$
\aligned
&z_1\mapsto e^{i\beta}z_1,\\
&z_j\mapsto e^{i\frac{\gamma}{\alpha_j}}z_j,\qquad j=2,\dots,n,
\endaligned
$$
where $\beta\in{\Bbb R}$, $-\pi<\gamma\le\pi$. Therefore,
for any such domain $D$, one has $\text{dim}\,\text{Aut}(D)\ge 4$.

If $n=2$, by differentiating both parts of (16) with respect to $t$ 
and $\overline{t}$ and setting $t=1$, one obtains that
$\psi(z_2)=c|z_2|^{\alpha}$, with $c>0$. Therefore, for $n=2$, domain (15)
is equivalent to a domain of the form (11) which is Reinhardt. However, as
examples in \cite{FIK2} show, there exist bounded domains in ${\Bbb C}^2$
with $C^{1,\beta}$-smooth boundary, for some $0<\beta<1$, with
non-compact automorphism group, that are not biholomorphically equivalent
to any Reinhardt domain and thus to any domain of the form (15). It would
be interesting to know if, for $k\ge 2$, there also exist $C^k$-smooth
bounded domains with non-compact automorphism group that are not
equivalent to any domain (15), or, for $k\ge 2$, that the domains (15) are,
in
fact, the only possibilities up to biholomorphic equivalence.

For comparison, we state below the conjecture of Bedford/Pinchuk
\cite{BP2} (see also \cite{BP1}) for domains with $C^\infty$-smooth
boundary. Assign weights $\alpha_j=2m_j$, $m_j\in{\Bbb N}$,
$j=2,\dots,n$, to the variables $\tilde z=(z_2,\dots,z_n)$. If
$K=(k_2,\dots,k_n)$ is a multi-index, we set
$\text{wt}(K)=\sum_{j=2}^n\frac{k_j}{\alpha_j}$. Consider real
polynomials of the form $$ P(\tilde z,\overline{\tilde
z})=\sum_{\text{wt}(K)=\text{wt}(L)=\frac{1}{2}}a_{KL}z^K\overline{\tilde
z}^L,\tag{17} $$ where $a_{KL}\in{\Bbb C}$ and
$a_{KL}=\overline{a_{LK}}$.  \medskip

{\bf Conjecture (Bedford/Pinchuk).} Any bounded domain with
non-compact automorphism group and $C^\infty$-smooth boundary is
biholomorphically equivalent to a domain $$ \left\{|z_1|^2+P(\tilde
z,\overline{\tilde z})<1\right\}, $$
where $P$ is a polynomial of the form (17).
\medskip

The above conjecture was proved in \cite{BP2} for convex domains of 
finite type and in \cite{BP1} for pseudoconvex domains of finite type for
which the Levi form of the boundary has rank at least $n-2$. Note that,
for polynomials (17), as well as for functions $\psi$ as in (5), condition
(16) is satisfied.

This work was initiated while the first 
author was an Alexander von Humboldt Fellow at the 
University of Wuppertal.  Research at MSRI by the second author was
supported in part by NSF Grant DMS-9022140 and also by
NSF Grant DMS-9531967.

\bigskip

{\obeylines
Centre for Mathematics and Its Applications 
The Australian National University 
Canberra, ACT 0200
AUSTRALIA 
E-mail address: Alexander.Isaev\@anu.edu.au
\smallskip
and
\smallskip
Bergische Universit\"at
Gesamthochschule Wuppertal
Mathematik (FB 07)
Gaussstrasse 20
42097 Wuppertal
GERMANY
E-mail address: Alexander.Isaev\@math.uni-wuppertal.de
\bigskip

Department of Mathematics
Washington University, St.Louis, MO 63130
USA 
E-mail address: sk\@math.wustl.edu
\smallskip
and
\smallskip
MSRI
1000 Centennial Drive
Berkeley, CA 94720
USA
E-mail address: krantz\@msri.org}

\enddocument